\documentclass[12pt, a4paper]{article}

\usepackage{amsmath}
\usepackage{amssymb}
\usepackage{amsthm}
\usepackage{graphicx}

\newtheorem{theorem}{Theorem}
\newtheorem{lemma}{Lemma}
\newtheorem{proposition}{Proposition}
\newtheorem{definition}{Definition}
\newtheorem{corollary}{Corollary}

\begin{document}
\noindent {\bf A global and superlinearly convergent algorithm for nonlinear nondifferential convex programming problems with a generalized Armijo line-search}\\

\noindent {\it Jiapu Zhang}\\

\noindent {\bf Abstract} This paper presents a new generalized Armijo line-search method, and combines it with a $\phi$-regulation defined to obtain a new algorithm solving the very general nonlinear nondifferential convex programming. For the algorithm designed, the global convergence is proved and the algorithm has superlinear convergent rate under very weak conditons. This paper genralized the results of reference ``Fukushima M and Qi LQ (1996) {\it SIAM J Optim} 6: 1106-20".

\section{Introduction}
It is well-known that nondifferentiable optimization is a very active area 
in the field of optimization. The study for its rapidly convergent rate 
motivates many feasible and efficient methods to be produced, such as the Cutting 
Plane method, the Subgradient method, the Bundle method, and the Proximal Desent 
method, etc.. However, all these methods at most have linear convergent rate. 
A method for the nonsmooth optimization problems that have superlinear 
convergence was presented recently in [1]. In fact, it is a globlized 
approximate Newton method. Even though computational experiments with this 
algorithm and the additional semismoothness assumptions are needed, it is 
the repaidest implimentable method by now to our knowledge.\\

In this paper, we present a globally and two-step superlinearly 
convergent algorithm, which can be used for nonlinear nondifferentiable convex 
programming in its wide sense. We had known well that the Armijo line-search 
method is often used in both the theoretical analysis and practical 
application for algorithms because of its excellent features. Study papers 
[2], [3], etc., we extend Armijo line search to a generalized Armijo line 
search, which is expansively described. The Moreau-Yosida regularization [4] 
is also extented to be $\phi$-regularization in this paper. Thirdly, the 
algorithm interation can be started from any given initial point in $R^n$. All 
these bases pave a way to construct a general algorithm for solving nonlinear 
nondifferentiable convex programming.\\

Studying the papers written by Moreau [4], Rockafellar [5,6,7], 
Hiriart-Urruty and Lemarechal [8], Fukushima and Qi [1] et al., we define 
a $\phi$-regularization for the objective function of a programme. We get that 
Moreau-Yosida regularization is a special case of $\phi$-regularization. At 
the meantime, we feel that semismoothness assuptions in [1] can be erased 
completely in this paper. The proofs for many conclusions in [1] can also be 
changed by the ways in this paper. We organized this paper as follows. In 
Section 2, we give the generalized Armijo line search method, and then define 
$\phi$-regular function and $\phi$-regularization of the objective function; 
also discuss their properties. In Section 3, we present an approximate algorith
m for nonlinear nondifferentiable convex programming and prove its feasibility.
 At last, we show that the proposed algorithm has a two-step Q-superlinear 
convergent rate under very weaker conditions compared to those of [1].

\section{A generalized Armijo line-search method and $\phi$-regularization}
In this paper, we consider the following unconstrained optimization problem
\begin{eqnarray}
\min_{x \in R^n} f(x),
\end{eqnarray}
where $f: R^n \longrightarrow R^1$ be a possibly nondifferentiable convex 
function, which is sloved by means of iterative methods, 
$x^{k+1} =x^k +\tau d^k (k=1,2, \dots )$, in which $x^0$ is any given 
starting point in $R^n$, 
$\tau_k$ denotes the iterative step and $d^k$ be the iterative direction.\\

We denote the gradient $\nabla f$ of $f$ by $g$, and introduce the following 
definition firstly.
\begin{definition} 
$[3]$The mapping $\varphi : [0,+ \infty )\longrightarrow [0,+\infty )$ is called 
{\rm forcing function}, if for any nonnegtive sequential numbers $\{ t_i \} $ 
then we have
\begin{equation*}
\lim_{i \longrightarrow \infty } \varphi (t_i) =0 \Rightarrow 
\lim_{i \longrightarrow \infty } t_i =0.
\end{equation*}
\end{definition}

Now we extend Armijo line-search to a generalized Armijo line-search, which is 
described as {\bf GA line-search method:}\\

\noindent {\it Let $a>0, \delta \in (0,1), \rho \in (0,1), c_1 >0$ be given constants, $M$ be a given 
positive integer, $\varphi$ be a forcing function and $\psi :[0,+\infty )
\longrightarrow [0,1)$ be a nonegative decreasing function; at the same time, 
$\psi (t)>0, \forall t>0.$ Denote $\gamma^k =-\frac{g_k^{T} d^k}
{\| d^k \|}$ and $\gamma_k =\min \{ \delta, c_1 \varphi (\gamma^k )+
\psi (\| g_k \| )\}$, select the step length $\tau_k =\rho^{h_k} a$, 
where $h_k$ is the first nonnegative integer $h$ satisfying
\begin{equation*}
f(x^k +\rho^h a d^k )\leq \max_{0 \leq j \leq m(k)} \{ f(x^{k-j} ) \} +\gamma_k
\rho^h a g_k^{T} d^k,
\end{equation*}
here $m(0)=0$ and $0\leq m(k) \min \{ m(k-1)+1, M \}.$}\\

Using this new line-search method and the following $\phi$-regularization and 
its properties, we shall study problem (1).

\begin{definition}
Assume that the function $\phi (\cdot ,\cdot ): R^n \times R^n \longrightarrow 
R^1$ satisfies the following properties:\\
(i) it is continuously differentiable on $R^n \times R^n$;\\
(ii) $\phi (z,x) \geq 0, \phi (z,x)=0 \Leftrightarrow z=x;$\\
(iii) it is strongly convex on $R^n \times R^n$, i.e., there exists a positive 
constant $\beta$ such that
\begin{eqnarray*}
%\begin{split}
\phi (z',x') -\phi (z,x) 
&\geq \nabla_z \phi (z,x)^{T} (z'-z) +
\nabla_x \phi (z,x)^{T} (x'-x)\\ 
&\quad +\beta \| (z'-z)-(x'-x) \|^2, \forall z',z,x',x \in R^n;
%\end{split}
\end{eqnarray*}
(iv) $\nabla_z \phi (z,x) $ and $\nabla_x \phi (z,x)$ are Lipschitz continuous 
on $R^n \times R^n$, i.e., there exists a positive real number $L$ such that
\begin{eqnarray*}
\| \nabla \phi (z',x') -\nabla \phi (z,x) \| \leq L \| (z'-z)-(x'-x) \|, \forall z',z,x',x \in R^n, 
\end{eqnarray*}
where $\nabla \phi (z,x)=(\nabla_z \phi (z,x),
\nabla_x \phi (z,x))^{T}.$\\
(v) $\nabla_z \phi (z,x)=- \nabla_x \phi (z,x).$\\
Then $\phi (\cdot ,\cdot )$ is called {\rm $\phi$-regular function} and 
$f_{\phi} =\min_{z \in R^n} \{ f(z) +{1 \over \lambda} \phi (z,x) \} 
(\lambda >0)$ is called {\rm $\phi$-regularization} of $f$.
\end{definition}      

By the above definition and the convexity of $f$, we can easily get $f(z)+{1 \over 
\lambda} \phi (z,x)$ is strongly convex and its level set is bounded; so 
the minimum is attained uniquely for each $x \in R^n$. We denote the unique 
minimizer by $p(x)$, i.e.,
\begin{eqnarray*}
f(p(x))+{1 \over \lambda} \phi (p(x),x) =\min_{z \in R^n} \{ f(z) +
{1 \over \lambda} \phi (z,x) \} . 
\end{eqnarray*}

The results stated in the following propositions are fundamental and useful 
in the subsequent discussions.

\begin{lemma}
Suppose $f_{\phi} (z,x)$ be two-unitary continuous function on 
$R^n \times R^n$ and satisfies the following conditions:\\
1). For any $z$, $f_\phi (z,x)$ be differentiable in respect of $x$;\\
2). $B(x)=\{ z | f_\phi (x) =\min_{z \in R^n} f_\phi (z,x) \}$ ;\\
3). $\lim_{t \downarrow 0} \frac{f_\phi (z,x+ts)- f_\phi (z,x)}{t} = 
\nabla_x f_\phi (z,x)^{T} s$ be uniformly convergent in respect of $z$;\\
Then we have $\partial f_\phi (x)= \overline{ conv\{ \nabla_x f_\phi (z,x) | 
z \in B(x) \} }$.
\end{lemma}
\noindent {\bf Proof.} This is the natural generalization of [9,Theorem 4, P112]. 
$\hfill{\Box}$\\

\begin{proposition}
The function $f_\phi$ is finite-valued, convex and everywhere differentiable 
with gradient $\nabla f_\phi (x)={1 \over \lambda} \nabla_x \phi (p(x),x).$ 
Moreover, the gradient mapping $\nabla f_\phi : R^n \longrightarrow R^n$ 
is Lipschitz continuous, i.e., there exists $c>0$ such that
\begin{eqnarray*}
\| \nabla f_\phi (x) -\nabla f_\phi (x') \| \leq c \| x'-x \|, 
\forall x',x \in R^n.
\end{eqnarray*}
\end{proposition}
\noindent {\bf Proof.} By (ii) and $f_\phi =f(p(x))+{1 \over \lambda} \phi (p(x),x), 
f_\phi$ is finite-valued.\\

Let $x=tx_1 +(1-t)x_2, \forall t \in [0,1], \forall x_1, x_2 \in R^n$, by the
 convexity of $f$ and (iii), we obtain
\begin{eqnarray*}
%\begin{split}
t f_\phi (x_1 ) + (1-t) f_\phi (x_2 )
&=&t f(p(x_1 ))+(1-t) f(p(x_2 ))\\
&+&{t \over \lambda } \phi (p(x_1 ),x_1 ) \quad + 
\frac{1-t}{\lambda } \phi (p(x_2 ),x_2 )\\
&\geq& f(tp(x_1 )+(1-t)p(x_2 ))+ {1 \over \lambda} \phi 
(tp(x_1 )\\
&+&(1-t)p(x_2 ),t x_1 +(1-t) x_2 )\\
&\geq& f(p(x))+{1 \over \lambda } \phi (p(x),x)\\
&=& f_\phi (x).
%\end{split}
\end{eqnarray*}
Hence, $f_\phi (x)$ is convex.\\

Let $f_\phi (z,x)=f(z)+{1 \over \lambda} \phi (z,x)$. Because
\begin{equation*}
\frac{f_\phi (z,x+ts)- f_\phi (z,x)}{t} = \frac{1}{t\lambda } 
[\phi(z,x+ts)- \phi (z,x)] =\frac{1}{\lambda} \nabla_x \phi (z,x)^{T} s +
o(\| s \| )
\end{equation*} 
implies
\begin{equation*}
\begin{split}
| \frac{f_\phi (z,x+ts)- f_\phi (z,x)}{t} -\nabla_x f_\phi (z,x)^{T} | 
&= | \frac{f_\phi (z,x+ts)- f_\phi (z,x)}{t} -{1 \over \lambda} \nabla_x 
\phi (z,x)^{T} s | \\
&= o( \| s \| ),
\end{split}
\end{equation*}
we have
\begin{equation*}
\lim_{t \downarrow 0} \frac{f_\phi (z,x+ts)- f_\phi (z,x)}{t} =\nabla_x f_\phi 
(z,x)^{T} s,
\end{equation*}
which is uniformly convergent with respect to $z$. Clearly, $f_\phi (z,x)$ satisfies
 conditions 1) 2) 3) of the Lemma 1. Therefore, by Lemma 1, we can obtain 
\begin{equation*}
\begin{split}
\partial f_\phi (x) &= \overline{conv \{ \nabla_z f_\phi (z,x) | z \in 
B(x)= \{ z | f_\phi (x)= \min_{z \in R^n} \{ f(z) + {1 \over \lambda } 
\phi (z,x) \} \} \} }\\
&= \overline{conv \{ \nabla_z f_\phi (z,x) | z \in B(x)=\{ p(x) \} \} }\\
&= \overline{conv \{ {1 \over \lambda } \nabla_x \phi (z,x) | 
z \in B(x)=\{ p(x) \} \} }\\
&= \overline{conv \{ {1 \over \lambda } \nabla_x \phi (p(x),x) \} }\\
&= {1 \over \lambda } \nabla_x \phi (p(x),x).
\end{split}
\end{equation*}
This explains $f_\phi (x)$ is everywhere differentiable with the unique 
gradient $\nabla f_\phi (x)={1 \over \lambda } \nabla_x \phi (p(x), x).$

Since $p(x)= var \min_{z \in R^n } \{ f(z)+ \frac{1}{\lambda } \phi (z,x) \}$, 
there exists $g \in \partial f(p(x))$ such that $0=g+\frac{1}{\lambda } 
\nabla_z \phi (p(x),x)$, i.e., $\nabla f_\phi (x) =g$, then we have
\begin{equation}
\langle \nabla f_\phi (x) - \nabla f_\phi (x'), p(x)-p(x') \rangle \geq 0, 
\forall x,x' \in R^n.
\end{equation}
By (iii), we have
\begin{equation}
\begin{split}
&\langle \nabla f_\phi (x) - \nabla f_\phi (x'),x-x' \rangle 
+ 
\langle \nabla f_\phi (x) - \nabla f_\phi (x'), p(x')-p(x) \rangle\\  
&\geq \frac{2\beta }{\lambda } \| (p(x)-p(x')) - (x-x') \|^2 .
\end{split}
\end{equation}
By (iv), for all $x,x' \in R^n$, we have
\begin{equation}
2 \lambda^2 \| \nabla f_\phi (x) - \nabla f_\phi (x') \|^2 \leq 
L^2 \| (p(x)-p(x'))-(x-x') \|^2.
\end{equation}
Combining (2)-(4), we get
\begin{equation}
\| \nabla f_\phi (x) - \nabla f_\phi (x') \|^2 
\leq \frac{L^2}{4\beta \lambda} \langle \nabla f_\phi (x) - \nabla f_\phi (x'),
x-x' \rangle .
\end{equation}  
Let $c=\frac{L^2}{4\beta \lambda}$, then by (5) we can obviously get
\begin{equation*}
\| \nabla f_\phi (x) - \nabla f_\phi (x') \| \leq c \| x- x' \|
\end{equation*}
which states $\nabla f_\phi$ is Lipschitz continuous. $\hfill{\Box}$

\begin{proposition}
The following statements are equivalent:\\
(A1) $x$ minimizes $f$;\\
(A2) $x=p(x);$\\
(A3) $\nabla f_\phi (x)=0;$\\
(A4) $x$ minimizes $f_\phi (x)$;\\
(A5) $f(x)=f(p(x));$\\
(A6) $f(x)=f_\phi (x).$
\end{proposition}
\noindent {\bf Proof.} (A1)$\Rightarrow$(A2) By (ii) and 
\begin{eqnarray*}
f(p(x)) \geq f(x) =f(x)+{1 \over \lambda } \phi (x,x) 
\geq f(p(x))+{1 \over \lambda } \phi (p(x),x),
\end{eqnarray*}
we get $p(x)=x$.\\

(A2)$\Rightarrow$(A3) By (ii) and Proposition 1, we have 
$\nabla f_\phi (x)={1 \over \lambda} \nabla_x \phi (x,x)=0.$\\

(A3)$\Rightarrow$(A4) Since $f_\phi (x)$ is convex, we have
\begin{eqnarray*}
f_\phi (x) =f_\phi (x) + \nabla f_\phi (x) (\bar{x} -x) 
\leq f_\phi (\bar{x} ), \forall \bar{x}.
\end{eqnarray*}
So we can easily get the result.\\

(A4)$\Rightarrow$(A5) By the differentiability and convexity of $f_\phi (x)$, 
we get $\nabla f_\phi (x)=0$. Hence, by Proposition 1, 
$\nabla_x \phi (p(x),x)=0.$ Also since 
\begin{eqnarray*}
%\begin{split}
\phi (p(x),p(x))- \phi (p(x),x) &\geq & \nabla_x \phi (p(x),x)^{T} (p(x)-x),\\ 
\phi (p(x),p(x)) &=&0,\\
\phi (p(x),x) &\geq & 0,
%\end{split}
\end{eqnarray*}
we have $x=p(x)$. So that $f(x)=f(p(x))$.\\

(A5)$\Rightarrow$(A6) Since 
\begin{eqnarray*}
f(x)=f(x)+\frac{1}{\lambda} \phi (x,x) \geq f_\phi (x) 
=f(p(x))+\frac{1}{\lambda} \phi (p(x),x)=f(x)+\frac{1}{\lambda} \phi (p(x),x)
\end{eqnarray*}
and
\begin{eqnarray*}
\phi (p(x),x) \geq 0,
\end{eqnarray*}
we have $\phi (p(x),x) \geq 0$, so that $f(x)=f_\phi (x)$.\\

(A6)$\Rightarrow$(A1) Since
\begin{eqnarray*}
%\begin{split}
f(x)&=& f_\phi (x),\\
f(x)&=& f(x)+\frac{1}{\lambda} \phi (x,x),\\
f_\phi (x)&=& f(p(x))+\frac{1}{\lambda} \phi (p(x),x),
%\end{split}
\end{eqnarray*}
and the uniqueness of $p(x)$, we obtain $x=p(x)$; hence $\nabla f_\phi (x)=0$. 
Then by Proposition 1, we get
\begin{eqnarray*}
f(z)=f(z)+\frac{1}{\lambda} \phi (z,z) \geq f_\phi (z) 
\geq f_\phi (x) + \nabla f_\phi (x)^{T} (z-x)=f_\phi (x)=f(x)
\end{eqnarray*}
for all $z \in R^n$. So, $x$ minimizes $f$. $\hfill{\Box}$\\

Proposition 2 states that the study of problem (1) can be transfered 
to solve problem
\begin{equation}
\min_{x \in R^n} f_\phi (x).
\end{equation}
However, $f_\phi (x)= \min_{z \in R^n} \{ f(z) + \frac{1}{\lambda} \phi (z,x) 
\}$ is difficult or even impossible to find an exact solution $p(x)$ to 
express $f_\phi (x)$. In practice, the approximation of $p(x)$ which is 
denoted by $p^a (x,\varepsilon )$ can be found by some implementable 
algorithms [10,11,12]. We suppose that for all $x \in R^n$ and any 
$\varepsilon >0$, there exists $p^a (x,\varepsilon )$ such that 
\begin{equation}
f(p^a (x,\varepsilon )+\frac{1}{\lambda} \phi (p^a (x,\varepsilon ),x) 
\leq f_\phi (x) +\varepsilon.
\end{equation}
With $p^a (x,\varepsilon )$, define approximations to $f_\phi (x)$ and 
$\nabla f_\phi (x)$ by
\begin{equation}
f_\phi^a (x,\varepsilon )=f(p^a (x,\varepsilon ))+\frac{1}{\lambda} 
\phi (p^a (x, \varepsilon )
\end{equation}
and
\begin{equation}
g^a (x,\varepsilon )=\frac{1}{\lambda} \nabla_x \phi (p^a (x,\varepsilon ),x)
\end{equation}
respectively. Then by the strongly convexity of $\theta (z)=f_\phi (z,x) =f(z)+\frac{1}
{\lambda} \phi (z,x)$,  and choose special value together with Proposition 1 
yield the following Lemma 2, in which you can see the approximate degrees what 
we had stated.
\begin{lemma}
Let $p^a (x,\varepsilon )$ be a vecter satisfying (7), and 
$f^a_\phi (x,\varepsilon )$ and $g^a (x,\varepsilon )$ be given by (8) 
and (9) respectively. Then we have
\begin{gather}
f_\phi (x) \leq f^a_\phi (x,\varepsilon ) \leq f_\phi (x) +\varepsilon,\\
\| p^a (x,\varepsilon )- p(x) \| \leq 
\sqrt{\frac{\lambda \varepsilon}{\beta}},\\
\| g^a (x,\varepsilon )-g(x) \| \leq \sqrt{\frac{L^2 \varepsilon}{\beta \lambda
}}.
\end{gather}
\end{lemma}

\section{An algorithm and its convergence}

{\bf \underline{The algorithm:}}\\

\noindent {\it {\small
{\bf Step 0} Let $a>0, \delta >0, \rho \in(0,1), c_1 >0$ be given constants, 
$M$ be a given positive integer. $\varphi$ be a forcing function and 
$\psi :[0,+\infty ) \longrightarrow [0,1)$ be a nonnegative decreasing 
function and $\psi (t)>0, \forall t>0.$ Choose any vector $x^0$ in $R^n$, 
and give a $\varepsilon_0 >0.$ Set $k :=0$;\\
{\bf Step 1} Compute $g^a (x^k ,\varepsilon_k )$, if $g^a (x^k ,\varepsilon_k )
=0$, then stop; otherwise, pick a positive semidefinite symmetric matrix 
$V_k \in R^{n \times n}$ and a scalar $\alpha_k >0$. Compute $d^k =- (V_k +
\alpha_k I )^{-1} g^a (x^k ,\varepsilon_k ).$\\
{\bf Step 2} Compute $\gamma^k =- \frac{g^a (x^k ,\varepsilon_k )^{T} d^k }{\| 
d^k \| }$ and $\gamma_k =\min \{ \delta ,c_1 \varphi (\gamma^k )+ \psi (\| g^a 
(x^k ,\varepsilon_k ) \| ) \}$, choose a scalar $\varepsilon_{k+1}$ such that 
$\varepsilon_{k+1} ={1 \over 2} \varepsilon_k$. Let $m(0)=0$ and 
$0 \leq m(k) \leq \min \{ m(k-1)+1,M \}$. Compute $f^a_\phi (x^{k-j} ,
\varepsilon_{k-j} ).$ Select the step length $\tau_k =\rho^{h_k} a$, where 
$h_k$ is the first nonnegative integer $h$ satisfying
\begin{eqnarray*}
f^a_\phi (x^k + \rho^h a d^k ,\varepsilon_{k+1} ) 
\leq \max_{0 \leq j \leq m(k)} \{ f^a_\phi (x^{k-j} ,\varepsilon_{k-j}) \} 
+ \gamma_k \rho^h a g^a (x^k ,\varepsilon_k )^{T} d^k +\varepsilon_k,
\end{eqnarray*} 
where $f^a_\phi (x^k +\rho^h ad^k ,\varepsilon_{k+1} )$ be computed by formula 
(8). Set $x^{k+1} :=x^k +\tau_k d^k$, let $k :=k+1$, and return to Step 1.
} }\\

Combining [1] and [13], using Lemma 2, we can easily get the next proposition 
which ensures the feasibility of the above algorithm.
\begin{proposition}
For every $k$, there exists $\overline{\tau_k} >0$ such that
\begin{equation*}
f^a_\phi (x^k +\tau d^k ,\varepsilon_{k+1} ) \leq \max_{0 \leq j \leq m(k)}
\{ f^a_\phi (x^{k-j} ,\varepsilon_{k-j} ) \} + \gamma_k \tau g^a (x^k ,
\varepsilon_k )^{T} d^k +\varepsilon_k
\end{equation*}
for all $\tau \in (0,\overline{\tau_k}),$ where $m(0)=0$ and $0 \leq m(k) \leq 
\min \{ m(k-1)+1,M \} ,\quad k \geq 1.$
\end{proposition} 

Combining [13] and [1], using the GA rule and its reverse side formula, by 
Lemma 2, Proposition 2, we can get the following theorem which establishes 
global convergence of the algorithm.
\begin{theorem}
Assume that the objective function $f$ of problem (1) is bounded from below. 
Let $\{ \alpha \}$ be a bounded sequence of positive numbers and suppose 
that the eigenvalues of matrix sequence $\{ V_k + \alpha_k I \}$ be 
uniformly bounded. Then any accumulation point of $\{ x^k \}$ generated by 
algorithm is an optimal solution of problem (1).
\end{theorem} 

\section{Superlinear convergent rate of the algorithm}
\begin{definition}
\cite{s14}
$g: R^n \longrightarrow R^n$ be a mapping. Denote $\Omega_g = \{ x \in R^n | g$ 
{\it is differentiable at} $x \}$, $\partial_B g(x)= \{ A \in R^{n \times n} | 
A= \lim_{x^k \rightarrow x} \nabla g(x^k),x^k \in \Omega_g \}$, we say 
that $g$ is {\rm BD-regular} at $x$ if $g$ is Lipschitz continuous and all 
matrices $A \in \partial_B g(x)$ are nonsingular.
\end{definition}

By contradiction, since $\partial_B g(x)$ be compact, we can obtain the 
following proposition.
\begin{proposition}
For each $x \in R^n$, every $V \in \partial_B g(x)$ is a symmetric positive 
semidefinite matrix. Moreover, if $g$ is BD-regular at $x$, then there 
exists a constant $c>0$ and a neighborhood $N$ of $x$ such that for all 
$y \in N$
\begin{equation*}
\langle d,Vd \rangle \geq c \| d \|^2, \forall d \in R^n, 
\forall V \in \partial_B g(y).
\end{equation*}
\end{proposition} 
\begin{lemma}
Suppose $\bar{x}$ be optimal solution of problem (6), the gradient mapping 
$g$ of $f_\phi (x)$ be BD-regular at $\bar{x}$, then $\bar{x}$ is the unique 
optimal solution of (6).
\end{lemma}
\noindent {\bf Proof.} We, by the contradiction method, assume that the optimal solution of (6) 
is not unique, then there exists a optimal solution sequence $\{ x^k \} ,
x^k \neq \bar{x}$ and $ x^k \longrightarrow \bar{x} (k \longrightarrow \infty )
, f_\phi (x^k)= f_\phi (\bar{x}).$\\
Since the Lipschitz function $g$ is differentiable almost everywhere, there 
exists $\{ z_k^{t_k} \subset \Omega_g$. For every $k$ fixed temporarily, we 
can choose $t_k$ sufficiently large enough such that
\begin{equation*}
\| z_k^{t_k} -x^k \| \leq \| x^k -\bar{x} \|^2.
\end{equation*}
Thus,
\begin{equation*}
\lim_{k \longrightarrow \infty} z_k^{t_k} = \bar{x}.
\end{equation*}
Noticing $f_\phi (\bar{x})=f_\phi (x^k), g(x^k)=0$ for any $k$, we have
\begin{equation}
\begin{split}
0&=f_\phi (\bar{x}) -f_\phi (x^k)\\
&=f_\phi (\bar{x}) -f_\phi (z_k^{t_k}) -( f_\phi (x^k)-f_\phi (z_k^{t_k}))\\
&=\langle g(z_-k^{t_k}),\bar{x} -z_k^{t_k} \rangle + \frac{1}{2} 
\langle \bar{x} -z_k^{t_k}, \nabla g(z_k^{t_k})^{T} (\bar{x}-z_k^{t_k}) \rangle
+o(\| \bar{x} -z_k^{t_k} \| )\\
&- \langle g(z_k^{t_k}), x^k-z_k^{t_k} \rangle -{1 \over 2}
\langle x^k -z_k^{t_k}, \nabla g(z_k^{t_k} )^{T} (x^k -z_k^{t_k} \rangle
+o(\| x_k -z_k^{t_k} \| ).
\end{split}
\end{equation}   
By
\begin{equation*}
\| z_k^{t_k} -\bar{x} \| \geq \| x^k -\bar{x} \| -\| x^k -z_k^{t_k} \| 
\geq \| x^k -\bar{x} \| ( 1-\| x^k-\bar{x} \| )
\end{equation*}
and for sufficiently large number $k$ 
\begin{equation*}
(1- \| x^k -\bar{x} \| )>0, 
\end{equation*} 
we get 
\begin{equation*}
\frac{\| x^k - z_k^{t_k} \| }{\bar{x} -z_k^{t_k} \| }
\leq \frac{\| x^k -\bar{x} \|^2 }{ \| x^k -\bar{x} \| (1-\| x^k -\bar{x} \| )}
\longrightarrow 0 (k \longrightarrow \infty ).
\end{equation*} 
Hence,
\begin{equation}
\| x^k -z_k^{t_k} \|^2 = o(\| \bar{x} -z_k^{t_k} \|^2 )
\end{equation}
\begin{equation*}
\begin{split}
| \langle g(z_k^{t_k} ),x^k - z_k^{t_k} \rangle |
&= | \langle g(z_k^{t_k} )-g(x^k),x^k-z_k^{t_k} \rangle | \\
&\leq \| g(z_k^{t_k} )-g(x^k) \| \cdot \|x^k -z_k^{t_k} \| \\
&\leq L \|x^k -z_k^{t_k} \|^2 ,
\end{split}
\end{equation*}
where $L$ be the Lipschitz constant.
So we have 
\begin{equation}
\langle g(z_k^{t_k} ), \bar{x} -z_k^{t_k} \rangle 
= o(\| \bar{x} -z_k^{t_k} \|^2 ).
\end{equation}
In the same way, we obtain
\begin{equation*}
\begin{split}
| \langle g(z_k^{t_k} ), \bar{x} -z_k^{t_k} \rangle | 
&= | \langle g(z_k^{t_k} )-g(\bar{x} ), \bar{x} -z_k^{t_k} \rangle | \\
&\leq \| g(z_k^{t_k} )-g(\bar{x}) \| \cdot \| \bar{x} -z_k^{t_k} \| \\
&\leq L \| \bar{x} -z_k^{t_k} \|^2 ,
\end{split}
\end{equation*}
and we also have
\begin{equation}
\langle g(z_k^{t_k} ), \bar{x} -z_k^{t_k} \rangle 
= o( \| \bar{x} -z_k^{t_k} \|^2 ),
\end{equation}
By (14)-(16) and
\begin{equation*}
| \langle x^k -z_k^{t_k} , \nabla g(z_k^{t_k} )^{T} (x^k -z_k^{t_k} ) \rangle |
 =o(\| \bar{x} -z_k^{t_k} \|^2 ),
\end{equation*} 
also by (13), we know
\begin{equation}
0=\frac{1}{2} \langle \bar{x} -z_k^{t_k} ,\nabla g(z_k^{t_k} )^{T} (\bar{x} -
z_k^{t_k} ) \rangle + ( \| \bar{x} - z_k^{t_k} \|^2 ).
\end{equation}
If we suppose 
\begin{equation*}
\begin{split}
&\lim_{k \longrightarrow \infty } \frac{\bar{x} -z_k^{t_k} }{\| \bar{x} -
z_k^{t_k} \| } =d \neq 0,\\
&\lim_{k \longrightarrow \infty } \nabla g(z_k^{t_k} )
=\bar{V} \in \partial_B g(\bar{x} ),
\end{split}
\end{equation*}
then by (17) we get 
\begin{equation*}
0=\frac{1}{2} 
\frac{
\langle \bar{x} -z_k^{t_k} , 
\nabla g(z_k^{t_k} )^{T} (\bar{x} -z_k^{t_k} ) \rangle
}
{\| \bar{x} -z_k^{t_k} \|^2} 
+\frac{o(\| \bar{x} -z_k^{t_k} \|^2 )}{\| \bar{x} -z_k^{t_k} \|^2 }.
\end{equation*}
Letting $k \longrightarrow \infty$, taking limit, we obtain 
\begin{equation*}
\langle d, \bar{V} d \rangle =0,
\end{equation*}
which is a contradiction for the positivity of $\bar{V}$. Thus, $\bar{x}$ is 
unique.  $\hfill{\Box}$ 
\begin{theorem}
Assume that the conditions of Theorem 1 be satisfied. If $\bar{x}$ be the 
optimal solution of problem (1) and $g$ is BD-regular at $\bar{x}$, 
then $\bar{x}$ is the unique optimal solution of (1) and that entire 
sequence $\{ x^k \}$ generated by algorithm converges to $\bar{x}$.
\end{theorem}
\noindent {\bf Proof.} By the convexity of $f_\phi (x)$ and Lemma 3, $\bar{x}$ is 
the unique optimal solution of problem (6). By Proposition 2, $\bar{x}$ is also 
the unique optimal solution of problem (1). Next, we prove that the entire 
sequence $\{ x^k \}$ converges to $\bar{x}$. In fact, because the optimal 
solution of $f_\phi (x)$ is unique we know that the level set 
$l=\{ x | f_\phi (x) \leq f_\phi (\bar{x} ) \} =\{ \bar{x} \}$
be bounded, i.e., the level set of $f_\phi (x) $ is bounded. So, for the 
iterative sequence $\{ x^k \}$ we have 
\begin{equation*}
f^a_\phi (x^k ,\varepsilon_k ) \leq \cdots \leq f^a_\phi (x^0 ,\varepsilon_0 )
+ \sum_i \varepsilon_i =q < +\infty.
\end{equation*}   
Combining with
\begin{equation*}
\sum_{k=1}^{\infty} \varepsilon_k <+\infty,
\end{equation*}
we get
\begin{equation*}
x^k \in l_q = \{ x | f_\phi (x) \leq q \}.
\end{equation*} 
Then because $l_q$ is bounded for any $q$, $\{ x^k \}$ be bounded
then $\{ x^k \}$ has no non-convergent subsequence, which implies 
$x^k \longrightarrow \bar{x} (k \longrightarrow \infty ).$ $\hfill{\Box}$\\ 

This sequentially convergent theorem is fundamental and useful in the 
subsequent discussions of the Q-superlinear convergence of the algorithm.\\

In the algorithm, alternatively, we may choose $V_k$ by calculating 
$\nabla g(\hat{x}^k )$, where $\hat{x}^k$ is very close to $x^k$ and $g$ 
is differentiable at $\hat{x}^k$. By the expression of set $\partial_B g(x), 
\nabla g(\hat{x}^k )$ can be made as close as possible to a member of 
$\partial_B g(x^k )$ by making $\hat{x}^k$ close to $x^k$.
We may also let $\alpha_k$ tend to zero as $k$ tends to infinity if we choose 
$V_k$ as an approximate member of $\partial_B g(x^k )$ because of 
Proposition 4. We also suppose $\nabla g$ is Lipschitz continuous in the 
neighbourhood $N$ of $x^k$[16]. Summarizing all the contents above, we 
establish the following two-times Q-superlinear convergence of the algorithm.
\begin{theorem}
Assume that the conditions of Theorem 2 be satisfied. Suppose 
furthermore that\\
1'. $\varepsilon_k =o(\| g(x^k )\|^2 )$;\\
2'.$\lim_{k \longrightarrow \infty} dist (V_k , \partial_B (x^k))=0,
\lim_{k \longrightarrow \infty} \alpha_k =0, and \quad \| V_k -\bar{V_k} \|
\leq l \| \hat{x}^k -\bar{x}^k \| ,where \quad V_k ,\bar{V_k} \in \partial_B
g(x^k ), \hat{x}^k \longrightarrow x^k ,\bar{x}^k \longrightarrow x^k;$\\
3'. $\tau_k \equiv 1$ for all large $k$, and $a=1$.\\
Then $\{ x^k \}$ converges to $\bar{x}$ two-step Q-superlinearly.
\end{theorem}
\noindent {\bf Proof.} First note that, by Theorem 2, the sequence $\{ x^k \}$ converges 
to $\bar{x}$. Then by the condition 1', the inequality (12), Proposition 1 
and 2, we get 
\begin{equation}
\begin{split}
\| g^a (x^k ,\varepsilon_k ) -g(x^k ) \|
&=O(\sqrt{\varepsilon_k })\\
&=o(\| g(x^k )\| )\\ 
&=o(\| g(x^k )-g(\bar{x} ) \| )\\
&=o(\| x^k -\bar{x} \| ).
\end{split}
\end{equation}
By the condition 2', there exists a $\bar{V_k} \in \partial_B g(x^k )$ such 
that
\begin{equation}
\| V_k -\bar{V_k} \| =o(1).
\end{equation}

We also have 
\begin{equation}
\| g(x^k )-g(\bar{x} ) -\bar{V_k } (x^k -\bar{x} )\| =o(\| x^k -\bar{x} \| )
\end{equation}
without the condition of $g$ being semismooth at $\bar{x}$. In fact, for each 
$k$ there exists $z_k^t \in \Omega_g$ satisfying 
\begin{equation*}
\lim_{t \longrightarrow \infty } z^t_k =x^k ,
\nabla g(z_k^t )\longrightarrow \bar{V_k } (t \longrightarrow \infty ),
\end{equation*}
and we can choose enough large $t_k$ such that
\begin{equation}
\| z_k^{t_k} -x^k \| \leq \| x^k -\bar{x} \|^2
\end{equation}
and
\begin{equation}
\| \nabla g(z_k^{t_k} )- \bar{V_k } \| \leq \| x^k -\bar{x} \| .
\end{equation}
Obviously,
\begin{equation*}
\lim_{k \longrightarrow \infty } z_k^{t_k} =\bar{x},
\end{equation*}
thus we have
\begin{equation}
\begin{split}
&g(x^k)-g(\bar{x})-\bar{V_k} (x^k -\bar{x})\\
&=g(x^k)-g(z_k^{t_k}) -(g(\bar{x})-g(z_k^{t_k}))-\bar{V_k}(x^k-\bar{x} )\\
&=g(x^k)-g(z_k^{t_k})+ \nabla g(z_k^{t_k})^{T} (z_k^{t_k} -\bar{x} )
-\bar{V_k} (x^k -\bar{x} )+o(\| z_k^{t_k} -\bar{x} \| )\\
&=g(x^k)-g(z_k^{t_k})+ \nabla g(z_k^{t_k}) (z_k^{t_k} -x^k )\\
&\quad +
(\nabla g(z_k^{t_k})-\bar{V_k})(x^k -\bar{x} )+o(\| z_k^{t_k} -\bar{x} \| ).
\end{split}
\end{equation} 
Since
\begin{equation*}
\| g(x^k)-g(z_k^{t_k} )\| =o(\| x^k -\bar{x} \| ),
\end{equation*}
we have 
\begin{equation}
\| g(x^k)-g(z_k^{t_k} ) \| =o(\| x^k -\bar{x} \| ).
\end{equation}
Because $g$ is BD-regularization at $\bar{x}$, $\lim_{k \longrightarrow \infty}
z_k^{t_k} =\bar{x}$, there exists $H>0$ such that
\begin{equation*}
\| \nabla g(z_k^{t_k}) \| \leq H.
\end{equation*}
Now by (21) and (22) we get
\begin{equation}
\| \nabla g(z_k^{t_k} )(x^k -z_k^{t_k} ) \| 
\leq H \| x^k - \bar{x} \|^2 , 
\end{equation}
\begin{equation}
\begin{split}
\| \nabla g(z_k^{t_k} )-(x^k - \bar{x} ) \| 
&\leq \| \nabla g(z_k^{t_k} )-
\bar{V_k} \| \cdot \| x^k -\bar{x} \| \\
&\leq \| x^k - \bar{x} \| \cdot 
\| x^k - \bar{x} \| = \| x^k - \bar{x} \|^2 ,
\end{split}
\end{equation}
\begin{equation}
\| z_k^{t_k} -\bar{x} \| 
\leq \| z_k^{t_k} - x^k \| + \| x^k - \bar{x} \| .
\end{equation} 
By (23)-(27), we can obtain (20) immediately. 
  
Notice  that $\| (V_k +\alpha_k I )^{-1} \| =O(1)$. Then by the algorithm, 
(18)-(20), the conditions $\lim_{k \longrightarrow \infty } \alpha_k =0$ and 
3', we have
\begin{equation*}
\begin{split}
\| x^{k+1} -\bar{x} \| 
&= \| x^k -\bar{x} -(V_k +\alpha_k I )^{-1} g^a (x^k ,\varepsilon_k ) \| \\
& \leq \| (V_k +\alpha_k I )^{-1} \| \cdot \| (V_k +\alpha_k I ) (x^k -
\bar{x} )-g^a (x^k ,\varepsilon_k ) \| \\
& \leq \| (V_k + \alpha_k I )^{-1} \| \{ \| g^a (x^k ,\varepsilon_k )-g(x^k ) 
\| \\
& \quad + \| g(x^k) -g(\bar{x} ) -\bar{V_k} (x^k -\bar{x} ) \| +
( \| V_k - \bar{V_k} \| + \alpha_k ) \| x^k -\bar{x} \| \} \\
&= o(\| x^k -\bar{x} \| ).
\end{split}
\end{equation*}  
for all large $k$. Therefore, 
\begin{equation*}
\lim_{k \longrightarrow \infty } \frac{\| x^{k+1} -\bar{x} \|}{\| x^k -\bar{x} 
\| } =0,
\end{equation*}
which means that $\{ x^k \}$ converges to $\bar{x}$ Q-superlinearly.

Now the two-step convergent rate is only left us to to prove. In fact, 
because of the BD-regularization of $g$ at $\bar{x}$, by Proposition 4, 
there exists $\hat{m} >0$ such that $\| x^k -\bar{x} \| \leq \hat{m}$ and 
there exists $r>0$ such that $\| V_k \| \leq r$ as $V_k$ being positive.
Since $\lim_{k \longrightarrow \infty} \alpha_k =0,$ we can let 
$\| (V_k +\alpha_k I )^{-1} \| \leq r'$. Hence, we have
\begin{equation*}
\begin{split}
x^{k+1} -\bar{x}
&= x^k +d^k -\bar{x} \\
&= x^k -(V_k +\alpha_k I )^{-1} g^a (x^k ,\varepsilon_k ) -\bar{x} \\
&= x^k -(V_k +\alpha_k I )^{-1} (g^a (x^k ,\varepsilon_k )- g(\bar{x} ))-
\bar{x} \\
&= x^k -\bar{x} -(V_k +\alpha_k I )^{-1} ( g^a (x^k ,\varepsilon_k )-
g(\bar{x}) ) \\
&= (V_k +\alpha_k I )^{-1} (V_k + \alpha_k I ) (x^k - \bar{x})-
(V_k + \alpha_k I )^{-1} ( g^a (x^k ,\varepsilon_k )- g(\bar{x} ) )\\
&= (V_k + \alpha_k I)^{-1} [ (V_k + \alpha_k I) (x^k -\bar{x} )
- ( g^a (x^k ,\varepsilon_k ) - g(\bar{x} ) )] \\
&= (V_k +\alpha_k I)^{-1} [(V_k + \alpha_k I )(x^k -\bar{x} )
-(\bar{V_k} +\alpha_k I )(x^k -\bar{x} ) ]\\
&= (V_k +\alpha_k I)^{-1} (V_k -\bar{V_k} ) (x^k -\bar{x} ).
\end{split}
\end{equation*} 
So, 
\begin{equation}
\begin{split}
\| x^{k+1} -\bar{x} \| 
&\leq \| (V_k +\alpha_k I )^{-1} \| \cdot l\theta \cdot
\| x^k -\bar{x} \| \cdot \| x^k -\bar{x} \| \\
&\leq r'l\theta \| x^k -\bar{x} \|^2  ,
\end{split}
\end{equation}   
where $V_k =\lim_{\hat{x}^k \longrightarrow x^k } \nabla g(\hat{x}^k ), 
\bar{V_k} = \lim_{\hat{x}^k \longrightarrow x^k } \nabla g (\hat{\hat{x}}^k),
\hat{\hat{x}}^k =\hat{x}^k +\theta (\hat{x}^k -\bar{x} ), \theta \in (0,1).$
By (28) and  $r'l\theta$ be a positive constant we get the conclusion the 
two-step convergent rate as $\{ x^k \}$ converges to $\bar{x}$. $\hfill{\Box}$\\

\begin{corollary}
Suppose that the conditions of Theorem 3 hold, except that the condition 1' is 
replaced by 1'': $\varepsilon_k =o(\| g( x^{k-1} ) \|^2 ).$ Then $\{ x^k \}$ 
converges to $\bar{x}$ at least 2-step Q-superlinearly.
\end{corollary}
\noindent {\bf Proof.} It is only to modify the proof of theorem 3 gentally, we can 
obtain this result. $\hfill{\Box}$\\

In fact, replacing $\delta <\frac{1}{2}$ by $\gamma_k <\frac{1}{2}$, 
proving similarlly as [1,Theorem 4], we also have the theorem below.
\begin{theorem}
Suppose that the conditions of Theorem 3, except the condition 3', hold
and 
$\delta <{1\over 2}, a=1$. Then the condition 3' holds and $\{ x^k \}$ 
converges to the unique solution $\bar{x}$ to problem (1) two-step Q-superlinearly. 
\end{theorem}

\section*{Acknowledgments}
{\small The author thanks Professor Wang CY for the discussions in his seminars}.

\end{document}